\documentclass[a4paper,english,oneside]{article}
\usepackage[T1]{fontenc}
\usepackage[utf8]{inputenc} 

\usepackage{geometry}
\usepackage{amsmath}
\usepackage{amsfonts}
\usepackage{amssymb}
\usepackage{amsthm}
\usepackage[english]{babel}
\usepackage{color}
\usepackage{lmodern}
\usepackage{pstricks-add}
\usepackage{hyperref}
\hypersetup{hidelinks}

\usepackage{mathscinet}
\usepackage{enumitem}
\usepackage[cal=boondoxo,bb=ams]{mathalfa}
\usepackage{microtype}

\usepackage{mathtools}
\usepackage{esint}

\theoremstyle{plain}
\newtheorem{theorem}{Theorem}[section]
\newtheorem{lemma}[theorem]{Lemma}
\newtheorem{corollary}[theorem]{Corollary}
\newtheorem{proposition}[theorem]{Proposition}

\theoremstyle{definition}
\newtheorem{definition}[theorem]{Definition}

\theoremstyle{remark}
\newtheorem{remark}{Remark}

    \DeclareMathOperator\tr{Tr}
     
\begin{document}

\title{On the blow -- up of solutions to semilinear damped wave equations with power nonlinearity in compact Lie groups}

\author{ Alessandro Palmieri} 

\date{}

\maketitle

\begin{abstract}

In this note, we prove a blow-up result for the semilinear damped wave equation in a compact Lie group with power nonlinearity $|u|^p$ for any $p>1$, under suitable integral sign assumptions for the initial data, by using an iteration argument. A byproduct of this method is the upper bound estimate for the lifespan of a local in time solution. As a preliminary result, a local (in time) existence result is proved in the energy space via Fourier analysis on compact Lie groups.

\end{abstract}

\begin{flushleft}
\end{flushleft}

\begin{flushleft}
\textbf{AMS Classification (2020)} Primary:  35B44, 35L71; Secondary: 43A30, 43A77, 58J45
\end{flushleft}

\section{Introduction}

Let $\mathbb{G}$ be a compact Lie group and let $\mathcal{L}$ be the Laplace -- Beltrami operator on $\mathbb{G}$ (which coincides with the Casimir element of the enveloping algebra). 
In the present work, we prove a blow-up result for the Cauchy problem for the semilinear damped wave equation with power nonlinearity, namely,
\begin{align}\label{semilinear CP damped}
\begin{cases} \partial_t^2 u-\mathcal{L} u+\partial_t u =| u|^p, &  x\in \mathbb{G}, \ t>0,\\
u(0,x)=\varepsilon u_0(x), & x\in \mathbb{G}, \\ \partial_t u(0,x)=\varepsilon u_1(x), & x\in \mathbb{G},
\end{cases}
\end{align} where $p>1$ and $\varepsilon$ is a positive constant describing the smallness of Cauchy data.

For the classical semilinear damped wave equation in $\mathbb{R}^n$ it has been proved in \cite{Mat76,TY01,Zhang01,IT05} that the critical exponent is the same one as for the semilinear heat equation, namely, the so -- called \emph{Fujita exponent} $p_{\mathrm{Fuj}}(n)\doteq 1+\frac{2}{n}$. In the pioneering paper \cite{Mat76} for the low dimensional cases $n=1,2$ the global existence of small data solutions in the supercritical case is proved. On the other hand, in \cite{TY01} the global existence  is proved in the supercritical case $p>p_{\mathrm{Fuj}}(n)$ for any $n\geqslant 1$ by working with compactly supported initial data together with the blow -- up of local in time solution (under suitable sign assumptions) in the subcritical case $1<p<p_{\mathrm{Fuj}}(n)$. Then, in \cite{Zhang01} a blow -- up result is proved even in the critical case $p=p_{\mathrm{Fuj}}(n)$ and, finally,  in \cite{IT05}  the globally existence in the supercritical case is proved without requiring compactly supported initial data and for all spatial dimensions.

In the non -- Euclidean framework, the semilinear damped wave equation has been studied also in the Heisenberg group $\mathbf{H}_n$ in \cite{Pal19,GP19DW}. The critical exponent for the semilinear damped wave equation in this nilpotent Lie group is the Fujita exponent $p_{\mathrm{Fuj}}(\mathcal{Q})$, where $\mathcal{Q}=2n+2$ is the homogeneous dimension of $\mathbf{H}_n$. This result is consistent with the critical exponent for semilinear heat equation in stratified Lie groups (see \cite{Pas98,RY18,GP19,GP19Car,Oka19}), which admits $p_{\mathrm{Fuj}}(\mathcal{Q})$ as critical exponent ($\mathcal{Q}$ being the homogeneous dimension of the stratified group).

Recently, the Cauchy problem for the semilinear heat equation with power nonlinearity has been studied in the framework of connected unimodular Lie groups in \cite{RY18}. In particular, in the compact case it has been proved for any exponent $p>1$ the nonexistence of global in time distributional solution, under certain sign assumptions for the initial data (cf. \cite[Remark 1.6]{RY18}). Purpose of the present work is to show an analogous result for the Cauchy problem \eqref{semilinear CP damped}. We may interpret this fact as follows: for a compact Lie group the Haar measure of the ball $B_R(e)$ around the identity element $e$ with respect to the Riemannian distance behaves as a constant as $R\to \infty$, that is, $\mathbb{G}$ has polynomial volume growth of order 0 (equivalently, the global dimension $D=D(\mathbb{G})\in\mathbb{N}$ of $\mathbb{G}$ is 0); so, formally, the critical exponent is $$\displaystyle{\lim_{D\to 0} p_{\mathrm{Fuj}}(D)=\infty}$$ meaning that a blow -- up result holds for any $p>1$.

Furthermore, before proving the blow -- up result we shall  prove a local in time existence result in the classical energy space by using representation theory for compact Lie groups. In particular, Plancherel formula plays a crucial role for the proof of the local existence result, since it allows us to work with the group Fourier transform by duality.

\subsection{Main results}

In what follows $L^q(\mathbb{G})$ denotes the space of $q$ -- summable functions on $\mathbb{G}$ with respect to the normalized Haar measure for $1\leqslant q <\infty$ (respectively, essentially bounded for  $q=\infty$) and for $s>0$ and $q\in(1,\infty)$ the Sobolev space $H^{s,q}_\mathcal{L}(\mathbb{G)}$ is defined as $$H^{s,q}_\mathcal{L}(\mathbb{G)}\doteq \left\{ f\in L^q(\mathbb{G}): (-\mathcal{L})^{s/2}f\in L^q(\mathbb{G})\right\}$$ equipped with the norm
$$\|f\|_{H^{s,q}_\mathcal{L}(\mathbb{G)}}\doteq\| f\|_{L^q(\mathbb{G})}+\| (-\mathcal{L})^{s/2} f\|_{L^q(\mathbb{G})}.$$ As customary, we denote the Hilbert space $H^{s,2}_\mathcal{L}(\mathbb{G)}$ simply by $H^{s}_\mathcal{L}(\mathbb{G)}$.

Let us begin by stating the local existence result for the semilinear Cauchy problem \eqref{semilinear CP damped}.

\begin{theorem}  \label{Thm loc esistence}
Let $\mathbb{G}$ be a compact, connected Lie group and let $n$ be the topological dimension of $\mathbb{G}$. Let us assume $n\geqslant 3$. Let $(u_0,u_1)\in H_\mathcal{L}^1(\mathbb{G})\times L^2(\mathbb{G})$ and $p>1$ such that $p\leqslant \frac{n}{n-2}$. Then, there exists $T=T(\varepsilon)>0$ such that the Cauchy problem \eqref{semilinear CP damped} admits a uniquely determined mild solution $$u\in \mathcal{C}\left([0,T],H^1_\mathcal{L}(\mathbb{G})\right)\cap \mathcal{C}^1\left([0,T],L^2(\mathbb{G})\right).$$

\end{theorem}


\begin{remark} \label{Remark ub p} The upper bound assumption for the exponent $p$ in Theorem \ref{Thm loc esistence} is made in order to apply a Gagliardo -- Nirenberg type inequality  proved in \cite[Remark 1.7]{RY19}. Even the restriction $n\geqslant 3$ is made to fulfill the assumptions for the employment of such inequality.
\end{remark}

\begin{remark} The assumption $n\geqslant 3$ in the statement of Theorem \ref{Thm loc esistence} is technical,  as we observed in Remark \ref{Remark ub p}. It can be removed by looking for solutions in weaker spaces than the  one in the statement of Theorem \ref{Thm loc esistence}, namely, in $$\mathcal{C}\left([0,T],H^s_\mathcal{L}(\mathbb{G})\right)\cap \mathcal{C}^1\left([0,T],L^2(\mathbb{G})\right)$$ for some $s\in (0,1)$.
\end{remark}

\begin{remark} \label{Remark L^1 reg intro} In the Euclidean case and in the Heisenberg group, the trick to get a global existence result for a not empty range for $p$ was to require additional $L^1$ regularity for the Cauchy data. In this way, one could improve the decay rates in the estimates of the $L^2$ norm of the solution to the corresponding linear homogeneous problem and its first order derivatives. However, in the compact case $L^2(\mathbb{G})\subset L^1(\mathbb{G})$ and we will see that by working with $L^1(\mathbb{G})$ -- regularity for $u_0,u_1$ no additional decay rate can be gained for the $L^2(\mathbb{G})$ -- norm of the solution of the corresponding homogeneous Cauchy problem (cf. Section \ref{Subsubsection L^1-L^2 est}).
\end{remark}

Before stating the blow -- up result, let us introduce a suitable notion of energy solutions  for the semilinear Cauchy problem \eqref{semilinear CP damped}.

\begin{definition} \label{Definition energy sol}
Let $(u_0,u_1)\in H_\mathcal{L}^1(\mathbb{G})\times L^2(\mathbb{G})$. We say that $$u\in \mathcal{C}\left([0,T),H^1_\mathcal{L}(\mathbb{G})\right)\cap \mathcal{C}^1\left([0,T),L^2(\mathbb{G})\right)\cap L^p_{\mathrm{loc}}\left([0,T)\times\mathbb{G}\right)$$ is an \emph{energy solution} on $[0,T)$ to \eqref{semilinear CP damped} if $u$ fulfills the integral relation
\begin{align}
& \int_{\mathbb{G}} \partial_t u(t,x) \psi(t,x)\, \mathrm{d}x -\int_{\mathbb{G}}  u(t,x) \psi_s(t,x)\, \mathrm{d}x+\int_{\mathbb{G}}  u(t,x) \psi(t,x)\, \mathrm{d}x \notag \\
 & \quad -\varepsilon \int_{\mathbb{G}} u_1(x) \psi(0,x)\, \mathrm{d}x +\varepsilon \int_{\mathbb{G}}  u_0(x) \psi_s(0,x)\, \mathrm{d}x-\varepsilon \int_{\mathbb{G}}  u_0(x) \psi(0,x)\, \mathrm{d}x \notag \\
 & \quad + \int_0^t\int_\mathbb{G} u(s,x) \big(\psi_{ss}(s,x)-\mathcal{L}\psi(s,x)-\psi_s(s,x)\big) \mathrm{d}x \, \mathrm{d}s = \int_0^t\int_\mathbb{G} |u(s,x)|^p \psi(s,x)\,  \mathrm{d}x \, \mathrm{d}s \label{def energ sol int relation}
\end{align} for any $\psi\in\mathcal{C}^\infty_0([0,T)\times \mathbb{G})$ and any $t\in (0,T)$.
\end{definition}

\begin{theorem}  \label{Thm blow up}
Let $\mathbb{G}$ be a compact Lie group. Let $(u_0,u_1)\in H_\mathcal{L}^1(\mathbb{G})\times L^2(\mathbb{G})$ be nonnegative and nontrivial functions and let $p>1$. Let $u\in \mathcal{C}\left([0,T),H^1_\mathcal{L}(\mathbb{G})\right)\cap \mathcal{C}^1\left([0,T),L^2(\mathbb{G})\right)\cap L^p_{\mathrm{loc}}\left([0,T)\times\mathbb{G}\right)$ be an energy solution to \eqref{semilinear CP damped} according to Definition \ref{Definition energy sol} with lifespan $T=T(\varepsilon)$.  Then, there exists a positive constant $\varepsilon_0=\varepsilon_0(u_0,u_1,p)>0$ such that for any $\varepsilon\in (0,\varepsilon_0]$ the energy solution $u$ blows up in finite time. Moreover, the upper bound estimate for the lifespan 
\begin{equation}
T(\varepsilon)\leqslant C \varepsilon^{-(p-1)}
\end{equation} holds, where the constant $C>0$ is independent of $\varepsilon$.
\end{theorem}

\subsection*{Notations} Throughout the paper we will employ the following notations: $\mathcal{L}$ denotes the Laplace -- Beltrami operator on $\mathbb{G}$; $\tr(A)= \sum_{j=1}^d a_{jj}$ and $A^*=(\overline{a_{ji}})_{1\leqslant i,j\leqslant d}$ denote the trace  and the adjoint matrix of $A=(a_{ij})_{1\leqslant i,j\leqslant d} \in \mathbb{C}^{d\times d}$, respectively; $\mathrm{d}x$ stands for the normalized Haar measure on the compact group $\mathbb{G}$; finally, we write $f\lesssim g$ when there exists a positive constant $C$ such that $f\leqslant Cg$ and $f\approx g$ when $g\lesssim f \lesssim g$. 

\section{Local  existence} \label{Section local existence}

In this section Theorem \ref{Thm loc esistence} will be proved. Let us recall briefly the notion of mild solutions to \eqref{semilinear CP damped}. We apply Duhamel's principle in order to represent the solution to the linear inhomogeneous  problem
\begin{align}\label{linear CP damped}
\begin{cases} \partial_t^2 u-\mathcal{L} u+\partial_t u =F(t,x), &  x\in \mathbb{G}, \ t>0,\\
u(0,x)= u_0(x), & x\in \mathbb{G}, \\ \partial_t u(0,x)= u_1(x), & x\in \mathbb{G}.
\end{cases}
\end{align} So, if we denote by $E_0(t,x)$ and $E_1(t,x)$ the fundamental solutions to \eqref{linear CP damped} in the homogeneous case $F= 0$ with initial data $(u_0,u_1)=(\delta_0,0)$ and $(u_0,u_1)=(0,\delta_0)$, respectively, then, we may represent the solution to \eqref{linear CP damped} as follows:
\begin{align*}
u(t,x) = u_0(x)\ast_{(x)} E_0(t,x) + u_1(x)\ast_{(x)} E_1(t,x)+ \int_0^t F(s,x)\ast_{(x)} E_1(t-s,x) \, \mathrm{d}s.
\end{align*}
Notice that in the previous representation formula we applied the invariance by time translations for the differential operator $\partial_t^2-\mathcal{L}+\partial_t$ and the identity $L\big(v\ast_{(x)}E_1(t,\cdot)\big)=v\ast_{(x)}L(E_1(t,\cdot))$ for any left -- invariant differential operator $L$ on $\mathbb{G}.$

Therefore, we say that $u$ is a \emph{mild solution} to \eqref{semilinear CP damped} on $[0,T]$ if $u$ is a fixed point for the nonlinear integral operator
\begin{align*}
N: u\in X(T) \to N u(t,x) \doteq \varepsilon u_0(x)\ast_{(x)} E_0(t,x) +\varepsilon
 u_1(x)\ast_{(x)} E_1(t,x)+ \int_0^t |u(s,x)|^p\ast_{(x)} E_1(t-s,x) \, \mathrm{d}s
\end{align*} in the evolution space $X(T)\doteq\mathcal{C}\left([0,T],H^1_\mathcal{L}(\mathbb{G})\right)\cap \mathcal{C}^1\left([0,T],L^2(\mathbb{G})\right)$, equipped with the norm
$$\|u\|_{X(T)}\doteq \sup_{t\in[0,T]} \left(\|u(t,\cdot)\|_{L^2(\mathbb{G})}+\|(-\mathcal{L})^{1/2}u(t,\cdot)\|_{L^2(\mathbb{G})}+\|\partial_t u(t,\cdot)\|_{L^2(\mathbb{G})}\right).$$

In order to show that $N$ admits a uniquely determined fixed point for $T=T(\varepsilon)>0$ sufficiently small, we will employ Banach's fixed point theorem. Nevertheless, before studying the semilinear problem, we have to deal with the corresponding linear homogeneous problem. In particular, we shall determine $L^2(\mathbb{G})$ -- $L^2(\mathbb{G})$ estimates via the group Fourier transform with respect to the spatial variable. Once these estimates will have been established, we might prove the local (in time) existence result applying a Gagliardo -- Nirenberg type inequality derived recently in \cite{RY19} in the more general frame of connected Lie groups (cf. Lemma \ref{Lemma GN ineq}). The section is organized as follows: in Section \ref{Subsection GFT} we recall the main tools from Fourier Analysis on compact Lie groups which are useful for our approach; then, in Section \ref{Subsection L^2-L^2 est} we will derive the $L^2(\mathbb{G})$ -- $L^2(\mathbb{G})$ estimates for the solution of the corresponding homogeneous linear problem and its first order derivatives by using Plancherel identity in the framework of compact Lie group; finally, in Section \ref{Subsection Fixed Point} it will be shown that the operator $N$ as a uniquely determined fixed point for $T$ sufficiently small.

\subsection{Group Fourier transform} \label{Subsection GFT}

In this section, we recall some results on Fourier Analysis on compact Lie groups. For a detailed presentation of this topic we refer to the monograph \cite{RT10} and to \cite[Chapter 2]{FR16}. 

A \emph{continuous unitary representation} $\xi:\mathbb{G}\to\mathbb{C}^{d_\xi\times d_\xi}$ of dimension $d_\xi$ is a continuous group homomorphism from $\mathbb{G}$ to the group of unitary matrix $\mathrm{U}(d_\xi,\mathbb{C})$, that is, $\xi(xy)=\xi(x)\xi(y)$ and $\xi(x)^*=\xi(x)^{-1}$  for all $x,y\in\mathbb{G}$ and the elements $\xi_{ij}:\mathbb{G}\to \mathbb{C}$ of the matrix representation $\xi$ are continuous functions for all $i,j\in \{1,\ldots,d_\xi\}$. Two representations $\xi,\eta$ of $\mathbb{G}$ are said \emph{equivalent} if there exists an invertible intertwining operator $A$ such that $A \xi(x) = \eta(x) A$ for any $x\in \mathbb{G}$. A subspace $W\subset \mathbb{C}^{d_\xi}$ is said $\xi$ -- invariant if $\xi(x) \cdot W\subset W$ for any $x\in \mathbb{G}$. The representation $\xi$ is said \emph{irreducible} if the only $\xi$ -- invariant subspaces are the trivial ones.

The unitary dual of $\mathbb{G}$, denoted by $\widehat{\mathbb{G}}$, consists of the equivalence class $[\xi]$ of continuous irreducible unitary representation $\xi:\mathbb{G}\to\mathbb{C}^{d_\xi\times d_\xi}$. 

For a function $f\in L^1(\mathbb{G})$ its Fourier coefficients at $[\xi]\in\widehat{\mathbb{G}}$ is defined by
\begin{align*}
\widehat{f}(\xi)\doteq \int_{\mathbb{G}} f(x) \xi(x)^*  \mathrm{d}x \in \mathbb{C}^{d_\xi\times d_\xi},
\end{align*} where the integral is taken with respect to the Haar measure on $\mathbb{G}$.

For $f\in L^2(\mathbb{G})$ the Fourier series representation is given by
$$f(x) = \sum_{[\xi]\in\widehat{\mathbb{G}}} d_\xi \tr\big(\xi(x)\widehat{f}\ \big)$$
for a.e. $x\in \mathbb{G}$, where hereafter just one unitary matrix representation is picked  in the sum for each equivalence class $[\xi]$ in $\widehat{\mathbb{G}}$.
Moreover, for $f\in L^2(\mathbb{G})$ Plancherel formula takes the following form
\begin{align} \label{Plancherel formula G}
\| f\|^2_{L^2(\mathbb{G})} =  \sum_{[\xi]\in\widehat{\mathbb{G}}} d_\xi \big\| \widehat{f}(\xi)\big\|^2_{\mathrm{HS}} \doteq \big\| \widehat{f} \, \big \|^2_{\ell^2(\widehat{\mathbb{G}})} 
\end{align} which allows us to introduce the norm on the Hilbert space $\ell^2(\widehat{\mathbb{G}})$ (cf. \cite[Section 10.3.3]{RT10} and \cite[Section 2.1.4]{FR16}), where the Hilbert -- Schmidt norm of the matrix $\widehat{f}(\xi)$ is given by
\begin{align*}
 \big\| \widehat{f}(\xi)\big\|^2_{\mathrm{HS}}  \doteq  \tr\big( \widehat{f}(\xi) \widehat{f}(\xi)^*\big) = \sum_{i,j=1}^{d_\xi} \big| \widehat{f}(\xi)_{ij}\big|^2.
\end{align*}

Let us analyze the behavior of the group Fourier transform for the Laplace -- Beltrami operator $\mathcal{L}$. Given $[\xi]\in\widehat{\mathbb{G}}$, then, all $\xi_{ij}$ are eigenfunctions of $\mathcal{L}$ with the same not positive eigenvalue $-\lambda_\xi^2$, that is,
\begin{align*}
-\mathcal{L} \xi_{ij}(x) = \lambda^2_\xi \, \xi_{ij}(x) \quad\mbox{for any} \ x\in\mathbb{G} \ \mbox{and any} \ i,j\in \{1,\ldots,d_\xi\}.
\end{align*} In other words, the symbol of $\mathcal{L}$ is 
\begin{align} \label{symbol Laplace-Beltrami}
\sigma_\mathcal{L} (\xi)= - \lambda^2_\xi I_{d_\xi},
\end{align} where $I_{d_\xi}\in \mathbb{C}^{d_\xi\times d_\xi}$ denotes the identity matrix, which means $\widehat{\mathcal{L} f}(\xi)= \sigma_\mathcal{L} (\xi) \widehat{f}(\xi) =- \lambda^2_\xi \widehat{f}(\xi)$ for any $[\xi]\in\widehat{\mathbb{G}}$.

Finally, by means of Plancherel formula for $s>0$ we may write
\begin{align*}
\|f\|^2_{\dot{H}^s_\mathcal{L}(\mathbb{G})} =\| (-\mathcal{L})^{s/2}f\|^2_{L^2(\mathbb{G})} =  \sum_{[\xi]\in\widehat{\mathbb{G}}} d_\xi \lambda_\xi^{2s}\big\| \widehat{f}(\xi)\big\|^2_{\mathrm{HS}}.
\end{align*}

\subsection{$L^2(\mathbb{G})$ -- $L^2(\mathbb{G})$ estimates for the solution to the  homogeneous problem} \label{Subsection L^2-L^2 est}

In this section, we derive $L^2(\mathbb{G})$ -- $L^2(\mathbb{G})$ estimates for the solution to the homogeneous problem associated to \eqref{semilinear CP damped}. We follow the main ideas from \cite{GR15}, namely, the group Fourier transform with respect to the spatial variable $x$ is applied together with Plancherel identity in order to determine by duality an explicit estimate for the $L^2(\mathbb{G})$ norms of $u(t,\cdot)$, $(-\mathcal{L})^{1/2}u(t,\cdot)$ and $\partial_t u(t,\cdot)$, respectively.

Let $u$ be a solution to \eqref{linear CP damped} in the homogeneous case $F=0$. So, let $\widehat{u}(t,\xi) = (\widehat{u}(t,\xi)_{k\ell})_{1\leqslant k,\ell\leqslant d_\xi}\in \mathbb{C}^{d_\xi\times d_\xi}$, $[\xi]\in\widehat{\mathbb{G}}$ denote  the group Fourier transform of $u$ with respect to the $x$ -- variable. Therefore, $\widehat{u}(t,\xi)$ solves the Cauchy problem for the system of ODEs (with size depending on the representation $\xi$)
\begin{align*}
\begin{cases}
\partial_t^2 \widehat{u}(t,\xi) -\sigma_\mathcal{L}(\xi) \widehat{u}(t,\xi)+\partial_t \widehat{u}(t,\xi) =0, & t>0, \\
 \widehat{u}(0,\xi) = \widehat{u}_0(\xi), \\
 \partial_t \widehat{u}(0,\xi) =  \widehat{u}_1(\xi).
\end{cases}
\end{align*} Thanks to \eqref{symbol Laplace-Beltrami}, the previous system is decoupled in $d_\xi^2$ independent scalar ODEs, namely,
\begin{align}\label{scalar dec ODE}
\begin{cases}
\partial_t^2 \widehat{u}(t,\xi)_{k\ell} +\partial_t \widehat{u}(t,\xi)_{k\ell} + \lambda_\xi^2 \widehat{u}(t,\xi)_{k\ell} =0, & t>0, \\
 \widehat{u}(0,\xi)_{k\ell} = \widehat{u}_0(\xi)_{k\ell}, \\
 \partial_t \widehat{u}(0,\xi)_{k\ell} =  \widehat{u}_1(\xi)_{k\ell},
\end{cases}
\end{align} for any $k,\ell\in \{1,\ldots,d_\xi\}$.

Straightforward computations lead to the following representation formula for the solution to the linear homogeneous Cauchy problem \eqref{scalar dec ODE}
\begin{align}\label{representation u hat kl}
  \widehat{u}(t,\xi)_{k\ell} = \mathrm{e}^{-\frac{t}{2}} G_0(t,\xi)  \widehat{u}_0(\xi)_{k\ell} +\mathrm{e}^{-\frac{t}{2}} G_1(t,\xi) \Big( \widehat{u}_1(\xi)_{k\ell}+ \tfrac 12  \widehat{u}_0(\xi)_{k\ell}\Big),
\end{align} where
\begin{equation} \label{def G0 G1}
\begin{split}
 G_0(t,\xi)  & \doteq \begin{cases} \cosh\left( \tfrac 12 \sqrt{1-4 \lambda_\xi^2} \, t\right)  & \mbox{if} \ \lambda_\xi^2<\tfrac 14 , \\
 1  & \mbox{if} \ \lambda_\xi^2=\tfrac 14, \\
 \cos\left( \tfrac 12 \sqrt{4 \lambda_\xi^2-1} \, t\right)  & \mbox{if} \ \lambda_\xi^2>\tfrac 14,   \end{cases} \\
 G_1(t,\xi)  & \doteq \begin{cases} \dfrac{2\sinh\left( \tfrac 12 \sqrt{1-4 \lambda_\xi^2} \, t\right)}{\sqrt{1-4 \lambda_\xi^2}}  & \mbox{if} \ \lambda_\xi^2<\tfrac 14 , \\
 t  & \mbox{if} \ \lambda_\xi^2=\tfrac 14, \\
 \dfrac{2\sin\left( \tfrac 12 \sqrt{4 \lambda_\xi^2-1} \, t\right)}{\sqrt{4 \lambda_\xi^2-1}}   & \mbox{if} \ \lambda_\xi^2>\tfrac 14.   \end{cases}
 \end{split}
\end{equation}

\subsubsection{Estimate for $\|u(t)\|_{L^2(\mathbb{G})}$} \label{Subsubsection u(t) L2-L2 est}
 Let us introduce the following partition of the unitary dual
\begin{align*}
D_1 & \doteq \left\{[\xi]\in \widehat{\mathbb{G}}: \lambda_\xi <\tfrac 18\right \}, \qquad D_2  \doteq \left\{[\xi]\in \widehat{\mathbb{G}}: \lambda_\xi \geqslant\tfrac 18\right \}.
\end{align*} Note that the choice of $1/8$ as threshold  in the previous definitions is irrelevant, since our goal is to separate $0$ (which is, for example, an eigenvalue for the continuous irreducible unitary representation $1:x\in\mathbb{G}\to 1\in \mathbb{C}$) from the other eigenvalues. For $[\xi]\in D_2$, it holds $$|\widehat{u}(t,\xi)_{k\ell}| \lesssim \mathrm{e}^{-ct}\big(|\widehat{u}_0(\xi)_{k\ell}| +|\widehat{u}_1(\xi)_{k\ell}|\big) \qquad \mbox{for any} \ t\geqslant 0,$$ where $c>0$ is a suitable constant independent of $[\xi]$. On the other hand, for $[\xi]\in D_1$ we can only get the estimate 
\begin{align} \label{estimate u hat kl}|\widehat{u}(t,\xi)_{k\ell}| \lesssim |\widehat{u}_0(\xi)_{k\ell}| +|\widehat{u}_1(\xi)_{k\ell}| \qquad \mbox{for any} \ t\geqslant 0.
\end{align} In Section \ref{Subsubsection L^1-L^2 est}, we will show that,  even if we worked with $L^1(\mathbb{G})$ -- regularity for $u_0,u_1$, the previous estimate cannot substantially by improved.

So, by using Plancherel formula, we obtain
\begin{align}
\| u(t,\cdot)\|^2_{L^2(\mathbb{G)}} = \sum_{[\xi]\in\widehat{\mathbb{G}}} d_\xi \sum_{k,\ell=1}^{d_\xi} |\widehat{u}(t,\xi)_{k\ell}|^2 & \lesssim  \sum_{[\xi]\in\widehat{\mathbb{G}}} d_\xi \sum_{k,\ell=1}^{d_\xi} \left( |\widehat{u}_0(\xi)_{k\ell}|^2 +|\widehat{u}_1(\xi)_{k\ell}|^2 \right) \notag \\ &  = \| u_0\|^2_{L^2(\mathbb{G)}} +\| u_1\|^2_{L^2(\mathbb{G)}}. \label{proof L^2 est u(t)}
\end{align}

\subsubsection{Estimate for $\|(-\mathcal{L})^{1/2} u(t)\|_{L^2(\mathbb{G})}$}

By Plancherel formula, we get
\begin{align*}
\| (-\mathcal{L})^{1/2} u(t,\cdot)\|^2_{L^2(\mathbb{G})} = \sum_{[\xi]\in\widehat{\mathbb{G}}} d_\xi \|\sigma_{(-\mathcal{L})^{1/2} }(\xi) \widehat{u}(t,\xi)\|^2_{\mathrm{HS}} = \sum_{[\xi]\in\widehat{\mathbb{G}}} d_\xi \sum_{k,\ell=1}^{d_\xi} \lambda_\xi^2 |\widehat{u}(t,\xi)_{k\ell}|^2.
\end{align*} For $[\xi]\in D_1$ it holds
\begin{align*}
\lambda_\xi^2 |\widehat{u}(t,\xi)_{k\ell}|^2 \lesssim \lambda_\xi^2 \mathrm{e}^{-2\lambda^2_\xi t} \left(|\widehat{u}_0(\xi)_{k\ell}|^2 +|\widehat{u}_1(\xi)_{k\ell}|^2\right)  \lesssim (1+t)^{-1} \left(|\widehat{u}_0(\xi)_{k\ell}|^2 +|\widehat{u}_1(\xi)_{k\ell}|^2\right),
\end{align*}  whereas for $[\xi]\in D_2$ it holds
\begin{align*}
\lambda_\xi^2 |\widehat{u}(t,\xi)_{k\ell}|^2 \lesssim \mathrm{e}^{- c t} \left( \lambda_\xi^2|\widehat{u}_0(\xi)_{k\ell}|^2 +|\widehat{u}_1(\xi)_{k\ell}|^2\right),
\end{align*} for a suitable positive constant $c$. Therefore,
\begin{align*}
\| (-\mathcal{L})^{1/2} u(t,\cdot)\|^2_{L^2(\mathbb{G})} & = \sum_{[\xi]\in D_1} d_\xi \sum_{k,\ell=1}^{d_\xi} \lambda_\xi^2 |\widehat{u}(t,\xi)_{k\ell}|^2 +\sum_{[\xi]\in D_2} d_\xi \sum_{k,\ell=1}^{d_\xi} \lambda_\xi^2 |\widehat{u}(t,\xi)_{k\ell}|^2 \\
& \lesssim (1+t)^{-1} \sum_{[\xi]\in D_1} d_\xi \sum_{k,\ell=1}^{d_\xi} \left(|\widehat{u}_0(\xi)_{k\ell}|^2 +|\widehat{u}_1(\xi)_{k\ell}|^2\right) \\ & \qquad + \mathrm{e}^{- c t} \sum_{[\xi]\in D_2} d_\xi \sum_{k,\ell=1}^{d_\xi} \left( \lambda_\xi^2|\widehat{u}_0(\xi)_{k\ell}|^2 +|\widehat{u}_1(\xi)_{k\ell}|^2\right) \\ & \lesssim  (1+t)^{-1} \left( \|u_0\|^2_{H^1_\mathcal{L}(\mathbb{G})}+\|u_1\|^2_{L^2(\mathbb{G})}\right).
\end{align*} Note that in the previous estimate the terms in the sum such that $[\xi]\in D_2$ provide the regularity for the Cauchy data.

\subsubsection{Estimate for $\|\partial_t u(t)\|_{L^2(\mathbb{G})}$}

Straightforward computations show that for any $[\xi]\in\widehat{\mathbb{G}}$ and any $k,\ell\in\{1,\ldots,d_\xi\}$ the following representation holds
\begin{align*}
\partial_t \widehat{u}(t,\xi)_{k\ell} = -\mathrm{e}^{-\frac{t}{2}} G_1(t,\xi) \lambda_\xi^2 \, \widehat{u}_0(\xi)_{k\ell} +\mathrm{e}^{-\frac{t}{2}} \left(G_0(t,\xi)-\tfrac 12 G_1(t,\xi) \right)  \widehat{u}_1(\xi)_{k\ell},
\end{align*} where $G_0(t,\xi),G_1(t,\xi)$ are defined in \eqref{def G0 G1}. In particular, for $[\xi]\in\widehat{\mathbb{G}}$ such that  $\lambda_\xi^2<\frac 18$ we may estimate
\begin{align*}
G_0(t,\xi)-\tfrac 12 G_1(t,\xi) & = \tfrac 12 \left(1-\tfrac{1}{\sqrt{1-4\lambda_\xi^2}}\right) \mathrm{e}^{(-\frac{1}{2}+\frac{1}{2}\sqrt{1-4\lambda_\xi^2})t} +\tfrac 12 \left(1-\tfrac{1}{\sqrt{1-4\lambda_\xi^2}}\right) \mathrm{e}^{(-\frac{1}{2}-\frac{1}{2}\sqrt{1-4\lambda_\xi^2})t} \\ 
& \approx  -\tfrac{1}{\sqrt{1-4\lambda_\xi^2}} \, \lambda_\xi^2 \, \mathrm{e}^{(-\frac{1}{2}+\frac{1}{2}\sqrt{1-4\lambda_\xi^2})t} +\tfrac 12 \left(1-\tfrac{1}{\sqrt{1-4\lambda_\xi^2}}\right) \mathrm{e}^{(-\frac{1}{2}-\frac{1}{2}\sqrt{1-4\lambda_\xi^2})t}.
\end{align*}
Combining Plancherel formula
\begin{align*}
\| \partial_t u(t,\cdot)\|^2_{L^2(\mathbb{G)}} = \sum_{[\xi]\in\widehat{\mathbb{G}}} d_\xi \sum_{k,\ell=1}^{d_\xi} |\partial_t\widehat{u}(t,\xi)_{k\ell}|^2 
\end{align*}
with the estimate
\begin{align*}
|\partial_t \widehat{u}(t,\xi)_{k\ell}| & \lesssim \lambda_\xi^2 \mathrm{e}^{-\lambda_\xi^2 t} \left(|\widehat{u}_0(\xi)_{k\ell}| +|\widehat{u}_1(\xi)_{k\ell}|\right) + \mathrm{e}^{-t} \left(|\widehat{u}_0(\xi)_{k\ell}| +|\widehat{u}_1(\xi)_{k\ell}|\right) \\
& \lesssim (1+t)^{-1} \left(|\widehat{u}_0(\xi)_{k\ell}| +|\widehat{u}_1(\xi)_{k\ell}|\right) 
\end{align*} for $[\xi]\in D_1$ and  the estimate
\begin{align*}
|\partial_t \widehat{u}(t,\xi)_{k\ell}| & \lesssim \mathrm{e}^{-c t} \left( \lambda_\xi|\widehat{u}_0(\xi)_{k\ell}| +|\widehat{u}_1(\xi)_{k\ell}|\right)
\end{align*} for $[\xi]\in D_2$, where $c>0$ is a suitable constant, then, we obtain
\begin{align*}
\| \partial_t u(t,\cdot)\|^2_{L^2(\mathbb{G)}} & \lesssim (1+t)^{-2} \sum_{[\xi]\in\widehat{\mathbb{G}}} d_\xi \sum_{k,\ell=1}^{d_\xi} \left( \lambda_\xi^2|\widehat{u}_0(\xi)_{k\ell}| ^2+|\widehat{u}_1(\xi)_{k\ell}|^2\right) = (1+t)^{-2}  \left( \|u_0\|^2_{H^1_\mathcal{L}(\mathbb{G})}+\|u_1\|^2_{L^2(\mathbb{G})}\right).
\end{align*}

Summarizing, in this section we proved the following result.
\begin{proposition} \label{Prop L^2-L^2 estimates}
Let us assume $(u_0,u_1)\in H^1_\mathcal{L}(\mathbb{G})\times L^2(\mathbb{G})$ and let $u\in \mathcal{C}\big([0,\infty),H^1_\mathcal{L}(\mathbb{G})\big)\cap  \mathcal{C}^1\big([0,\infty),L^2(\mathbb{G})\big)$ be the solution to the homogeneous Cauchy problem
\begin{align} \label{linear CP damped hom}
\begin{cases} \partial_t^2 u-\mathcal{L} u+\partial_t u =0, &  x\in \mathbb{G}, \ t>0,\\
u(0,x)= u_0(x), & x\in \mathbb{G}, \\ \partial_t u(0,x)= u_1(x), & x\in \mathbb{G}.
\end{cases}
\end{align} Then, $u$ satisfies the following $L^2(\mathbb{G})$ -- $L^2(\mathbb{G})$ estimates
\begin{align}
\|u(t,\cdot)\|_{L^2(\mathbb{G})} &\leqslant C \left( \|u_0\|_{L^2(\mathbb{G})}+\|u_1\|_{L^2(\mathbb{G})}\right), \label{L^2 norm u(t)} \\
\|(-\mathcal{L})^{1/2}u(t,\cdot)\|_{L^2(\mathbb{G})} &\leqslant C (1+t)^{-\frac{1}{2}}\left( \|u_0\|_{H^1_\mathcal{L}(\mathbb{G})}+\|u_1\|_{L^2(\mathbb{G})}\right), \label{L^2 norm (-L)^1/2 u(t)} \\
\| \partial_t u(t,\cdot)\|_{L^2(\mathbb{G})} &\leqslant C (1+t)^{-1}\left( \|u_0\|_{H^1_\mathcal{L}(\mathbb{G})}+\|u_1\|_{L^2(\mathbb{G})}\right), \label{L^2 norm ut(t)}
\end{align} for any $t\geqslant 0$, where $C$ is a positive multiplicative constant.
\end{proposition}

\subsubsection{$L^1(\mathbb{G})$ -- $L^2(\mathbb{G})$ estimates for the solution to the  homogeneous problem} \label{Subsubsection L^1-L^2 est} 

In Proposition \ref{Prop L^2-L^2 estimates} we employed data on $L^2(\mathbb{G})$ basis. Nevertheless, the  continuous embedding $L^2(\mathbb{G})\hookrightarrow L^1(\mathbb{G})$ holds due to the fact that Haar measure on the compact Lie group $\mathbb{G}$ is finite, so we might wonder what happens if we employ $L^1(\mathbb{G})$ -- regularity for $u_0,u_1$ instead.

For this purpose, let us recall first the definition of the space $\ell^\infty(\widehat{\mathbb{G}})$. Let $\mathcal{S}'(\widehat{\mathbb{G}})$ denotes the space of slowly increasing distribution on the unitary dual, whose definition can be found in \cite[Section 10.3.2]{RT10} or in \cite[Section 2.1.3]{FR16}. 

The space $\ell^\infty(\widehat{\mathbb{G}})$ is the subspace of $\mathcal{S}'(\widehat{\mathbb{G}})$ consisting of the functions $$H=\big\{H([\xi])\big\}_{[\xi]\in\widehat{\mathbb{G}}}: \ H([\xi])\in\mathbb{C}^{d_\xi\times d_\xi} \quad \mbox{for any} \ [\xi]\in \widehat{\mathbb{G}}$$ such that
\begin{align*}
\| H\|_{\ell^\infty(\widehat{\mathbb{G}})}\doteq \sup_{[\xi]\in\widehat{\mathbb{G}}} d_\xi^{-\frac 12} \| H(\xi)\|_{\mathrm{HS}}<\infty.
\end{align*} We refer to \cite[Section 10.3.3]{RT10} or \cite[Section 2.1.4]{FR16} for further details on the construction of the space $\ell^\infty(\widehat{\mathbb{G}})$.

 For the group Fourier transform it holds
\begin{align} \label{Riemann Lebesgue ineq on G}
\| \widehat{f}\|_{\ell^\infty(\widehat{\mathbb{G}})} \leqslant \| f\|_{L^1(\mathbb{G})}
\end{align} for any $f\in L^1(\mathbb{G})$ (cf. \cite[Proposition 10.3.42]{RT10}). Let us stress that  \eqref{Riemann Lebesgue ineq on G} is crucial if we want to use the $L^1(\mathbb{G})$ -- regularity for the Cauchy data.

With the same notations as in Section \ref{Subsubsection u(t) L2-L2 est}, it holds
\begin{align*}
\sum_{[\xi]\in D_2} d_\xi \sum_{k,\ell=1}^{d_\xi} |\widehat{u}(t,\xi)_{k\ell}|^2 & \lesssim  \mathrm{e}^{-2ct} \sum_{[\xi]\in D_2} d_\xi \sum_{k,\ell=1}^{d_\xi}\big(|\widehat{u}_0(\xi)_{k\ell}|^2 +|\widehat{u}_1(\xi)_{k\ell}|^2\big) \\
 & \lesssim  \mathrm{e}^{-2ct} \left(\| u_0\|^2_{L^2(\mathbb{G)}} +\| u_1\|^2_{L^2(\mathbb{G)}}\right)
\end{align*} for a suitable positive constant $c$. Therefore,  the addends such that $[\xi]\in D_1$ in \eqref{proof L^2 est u(t)} are the ones that do not provide a decay rate for $\|u(t,\cdot)\|_{L^2(\mathbb{G})}$. So, if we want to use $L^1(\mathbb{G})$ -- regularity in place of $L^2(\mathbb{G})$ -- regularity, then, necessarily, we have to apply it in the estimate of the terms with $[\xi]\in D_1$.


 The best estimates that we can obtain for any $[\xi]\in D_1$ for the multiplier in \eqref{representation u hat kl}  are  $$\mathrm{e}^{-\frac{t}{2}} G_0(t,\xi),\mathrm{e}^{-\frac{t}{2}} G_1(t,\xi)\leqslant \mathrm{e}^{-\lambda_\xi^2 t}.$$ The previous inequalities imply in turn
 \begin{align*}
 \sum_{[\xi]\in D_1} d_\xi \sum_{k,\ell=1}^{d_\xi} |\widehat{u}(t,\xi)_{k\ell}|^2  
 & \lesssim  \sum_{[\xi]\in D_1 } d_\xi \, \mathrm{e}^{-2\lambda^2_\xi t} \sum_{k,\ell=1}^{d_\xi} \big(|\widehat{u}_0(\xi)_{k\ell}|^2 +|\widehat{u}_1(\xi)_{k\ell}|^2\big)  \\
 & \lesssim  \sum_{[\xi]\in D_1 } d_\xi \, \mathrm{e}^{-2\lambda^2_\xi t}  \big(\|\widehat{u}_0(\xi)\|^2_{\mathrm{HS}} +\|\widehat{u}_1(\xi)\|^2_{\mathrm{HS}}\big).
 \end{align*} By using \eqref{Riemann Lebesgue ineq on G}, we have
 \begin{align*}
 \sum_{[\xi]\in D_1} d_\xi \, \mathrm{e}^{-2\lambda^2_\xi t}  \big(\|\widehat{u}_0(\xi)\|^2_{\mathrm{HS}} +\|\widehat{u}_1(\xi)\|^2_{\mathrm{HS}}\big)  & \lesssim  \Big(\sup_{[\xi]\in \widehat{\mathbb{G}}} d_\xi^{-\frac{1}{2}} \big(\|\widehat{u}_0(\xi)\|_{\mathrm{HS}} +\|\widehat{u}_1(\xi)\|_{\mathrm{HS}}\big)\Big)^2   \sum_{[\xi]\in D_1 } d_\xi^2 \, \mathrm{e}^{-2\lambda^2_\xi t}  \\
  & \lesssim  \Big(\|\widehat{u}_0\|_{\ell^\infty(\widehat{\mathbb{G}})} +\|\widehat{u}_1\|_{\ell^\infty(\widehat{\mathbb{G}})} \Big)^2   \sum_{[\xi]\in D_1} d_\xi^2 \, \mathrm{e}^{-2\lambda^2_\xi t} \\
 &  \lesssim  \Big(\|u_0\|_{L^1(\mathbb{G})} +\|u_1\|_{L^1(\mathbb{G})}  \Big)^2   \sum_{[\xi]\in D_1 } d_\xi^2 \, \mathrm{e}^{-2\lambda^2_\xi t}.
 \end{align*} 
Since the spectrum of $-\mathcal{L}$ is discrete (with finite dimensional eigenspaces) and has no finite cluster point, the sum in the right -- hand side of the previous chain of inequalities is a finite sum. However, since we have at least one continuous irreducible unitary representation such that $\lambda^2_\xi=0$ (namely, the trivial representation $1:x\in\mathbb{G}\to 1\in \mathbb{C}$), this sum cannot provide any decay rate.

Summarizing, the reason why we cannot get any decay rate for the norm $\|u(t,\cdot)\|_{L^2(\mathbb{G})}$ is that the Plancherel measure in the compact case is a weighted counting measure which does not allow to neglect the eigenvalue $0$. 

\subsection{Proof of Theorem \ref{Thm loc esistence}} \label{Subsection Fixed Point}

A fundamental tool to prove the local existence result is the following Gagliardo -- Nirenberg type inequality, whose proof can be found in \cite{RY19}.

\begin{lemma}\label{Lemma GN ineq}
Let $\mathbb{G}$ be a connected unimodular Lie group with topological dimension $n$. For any $1<q_0<\infty$, $0<q,q_1<\infty$ and $0<\alpha<n$ such that $q_0<\frac{n}{\alpha}$ the following Gagliardo -- Nirenberg type inequality holds
\begin{align} \label{GN ineq}
\| f\|_{L^q(\mathbb{G})}\lesssim \| f\|^\theta_{H^{\alpha,q_0}_\mathcal{L}(\mathbb{G})} \| f\|^{1-\theta}_{L^{q_1}(\mathbb{G})}
\end{align}
for any $f\in H^{\alpha,q_0}_\mathcal{L}(\mathbb{G}) \cap L^{q_1}(\mathbb{G})$, provided that $$\theta=\theta(n,\alpha,q,q_0,q_1)\doteq \frac{\frac{1}{q_1}-\frac{1}{q}}{\frac{1}{q_1}-\frac{1}{q_0}+\frac{\alpha}{n}}\in [0,1].$$
\end{lemma}

\begin{remark} In \cite[Remark 1.7]{RY19} the inequality \eqref{GN ineq} is provided in the more general framework of connected Lie group employing an absolutely continuous  measure with respect to the Haar measure, which coincides with the Haar measure in the case of unimodular Lie groups. Moreover, in \cite{RY19} the authors work with Sobolev spaces defined via a system of left -- invariant vector fields fulfilling H\"ormander's bracket generating condition. Clearly, for the Laplace -- Beltrami operator these Sobolev spaces coincide with the Sobolev spaces we recall in the introduction. 
\end{remark}

\begin{remark} \label{Remark local dimension}
Let us remind briefly the notion of local dimension for $\mathbb{G}$ with respect a system of left -- invariant vector fields and explain why we may restrict ourselves to work with the topological dimension of $\mathbb{G}$ in the statement of Lemma \ref{Lemma GN ineq} differently from \cite{RY19}. Let $X=\{X_1,\ldots,X_k\}$ be a system of left -- invariant vector field on $\mathbb{G}$ satisfying H\"ormander's condition. Let us consider the corresponding sub -- Laplacian $\mathcal{L}_X= \sum_{j=1}^k X_j^2 $. The local dimension of $\mathbb{G}$ with respect to the system $X$ is the natural number $d=d(X)$ such that the Haar measure $V_R$ of the ball $B_R(e)$ with respect to the Carnot -- Carath\`{e}odory distance satisfies the estimate $$c_1 R^d\leqslant V_R\leqslant c_2 R^d \quad \mbox{for any} \ R\in (0,1) $$ for some positive constants $c_1,c_2$. For further details on the local dimension we refer to \cite{NSW85,Var88}. In our case, since we are working with the Laplace -- Beltrami operator, which can be written as the sum of squares for a basis of the Lie algebra, the local dimension is nothing but the topological dimension of $\mathbb{G}$ (see Section II.4 in \cite{DER03}).
\end{remark}

\begin{remark}
Note that we can include the case $\theta(n,\alpha,q,q_0,q_1)=0$ in the statement of Lemma \ref{Lemma GN ineq} (which is not included in \cite{RY19}), since this corresponds to the trivial case $q=q_1$.
\end{remark}

\begin{corollary}\label{Lemma GN ineq L2}
Let $\mathbb{G}$ be a connected unimodular Lie group with topological dimension $n\geqslant 3$. For any $q\geqslant 2$ such that $q\leqslant \frac{2n}{n-2}$ the following Gagliardo -- Nirenberg type inequality holds
\begin{align} \label{GN ineq L2}
\| f\|_{L^q(\mathbb{G})}\lesssim \| f\|^{\theta(n,q)}_{H^{1}_\mathcal{L}(\mathbb{G})} \| f\|^{1-\theta(n,q)}_{L^{2}(\mathbb{G})}
\end{align}
for any $f\in H^{1}_\mathcal{L}(\mathbb{G})$, where $ \theta(n,q)\doteq  n\left(\frac{1}{2}-\frac{1}{q}\right)$.
\end{corollary}

We can now prove the existence of a uniquely determined solution to \eqref{semilinear CP damped} in $X(T)$ for $T$ sufficiently small.

Let us estimate $\|Nu\|_{X(T)}$ for $u\in X(T)$. We begin by  rewriting $Nu=  u^{\mathrm{ln}}+ J u$, where $$ u^{\mathrm{ln}}(t,x)\doteq \varepsilon u_0(x)\ast_{(x)} E_0(t,x) +\varepsilon  u_1(x)\ast_{(x)} E_1(t,x) $$ and $$J u(t,x) \doteq \int_0^t |u(s,x)|^p\ast_{(x)} E_1(t-s,x) \, \mathrm{d}s.$$ By Proposition \ref{Prop L^2-L^2 estimates} it follows immediately $\|u^{\mathrm{ln}}\|_{X(T)}\lesssim  \varepsilon\,  \|(u_0,u_1)\|_{H^1_\mathcal{L}(\mathbb{G})\times L^2(\mathbb{G})}$. On the other hand, thanks to the invariance by time translations of the linear Cauchy problem \eqref{linear CP damped hom}, we get 
\begin{align}
\|\partial_t^j (-\mathcal{L})^{i/2} Ju(t,\cdot)\|_{L^2(\mathbb{G})}  & \lesssim \int_0^t (1+t-s)^{-j-\frac{i}{2}} \| u(s,\cdot)\|^p_{L^{2p}(\mathbb{G})} \, \mathrm{d}s \notag \\ 
& \lesssim \int_0^t \| u(s,\cdot)\|^{p\theta(n,2p)}_{H^{1}_{\mathcal{L}}(\mathbb{G})} \| u(s,\cdot)\|^{p(1-\theta(n,2p))}_{L^2(\mathbb{G})} \, \mathrm{d}s \lesssim t \, \| u\|_{X(t)}^p \label{estimate Ju in X(T)}
\end{align} for $i,j\in\{0,1\}$ such that $0\leqslant i+j\leqslant 1$. Notice that the employment of \eqref{GN ineq L2}  in the previous estimate is the reason why we required the upper bound for $p$ ($p\leqslant \frac{n}{n-2}$) in the statement of Theorem \ref{Thm loc esistence}. Analogously, combining H\"older's inequality and \eqref{GN ineq L2}, for $i,j\in\{0,1\}$ such that $0\leqslant i+j\leqslant 1$ we obtain
\begin{align}
\|\partial_t^j (-\mathcal{L})^{i/2} (Ju(t,\cdot)-Jv(t,\cdot))\|_{L^2(\mathbb{G})}  & \lesssim \int_0^t (1+t-s)^{-j-\frac{i}{2}} \| |u(s,\cdot)|^p-|v(s,\cdot)|^p\|_{L^{2}(\mathbb{G})} \, \mathrm{d}s \notag \\  
& \lesssim \int_0^t  \| u(s,\cdot)-v(s,\cdot)\|_{L^{2p}(\mathbb{G})} \left(\|u(s,\cdot)\|^{p-1}_{L^{2p}(\mathbb{G})}+ \|v(s,\cdot)\|^{p-1}_{L^{2p}(\mathbb{G})}\right) \, \mathrm{d}s \notag \\ 
& \lesssim t \, \| u-v\|_{X(t)}\left(\|u\|^{p-1}_{X(t)}+\|v\|^{p-1}_{X(t)}\right). \label{estimate Ju -Jv in X(T)}
\end{align} Summarizing, we proved that 
\begin{align*}
\| N u\|_{X(T)} &\leqslant C  \varepsilon\,  \|(u_0,u_1)\|_{H^1_\mathcal{L}(\mathbb{G})\times L^2(\mathbb{G})} +C T \|u\|^p_{X(T)}, \\
\| N u -Nv\|_{X(T)} &\leqslant C T \| u-v\|_{X(T)}\left(\|u\|^{p-1}_{X(T)}+\|v\|^{p-1}_{X(T)}\right).
\end{align*} Therefore, for $T$ sufficiently small $N$ is a contraction on a certain ball around $0$ in the Banach space $X(T)$, so Banach's fixed point provides a uniquely determined fixed point $u$ for $N$ which is exactly our mild solution to \eqref{semilinear CP damped} on $[0,T]$.

\section{Blow -- up result}

In this section we are going to prove Theorem \ref{Thm blow up} by using an iteration argument. In particular, we will use a slicing procedure that allows us to treat an unbounded exponential multiplier that was introduced for the first time in \cite{ChenPal19MGT} and then applied to several different semilinear hyperbolic models (see \cite{ChenPal19MGTder,ChenRei20,Chen20,Chen20N}).

Let $u$ be a local in time energy solution to \eqref{semilinear CP damped} according to Definition \ref{Definition energy sol} with lifespan $T$. Let us fix $t\in (0,T)$. We can choose a bump function $\psi\in \mathcal{C}^\infty_0 ([0,T)\times \mathbb{G})$ such that $\psi=1$ on $[0,t]\times\mathbb{G}$ in \eqref{def energ sol int relation}. Then,
\begin{align*}
& \int_{\mathbb{G}} \partial_t u(t,x) \, \mathrm{d}x +\int_{\mathbb{G}}  u(t,x) \, \mathrm{d}x  -\varepsilon \int_{\mathbb{G}} u_1(x) \, \mathrm{d}x -\varepsilon \int_{\mathbb{G}}  u_0(x) \, \mathrm{d}x  = \int_0^t\int_\mathbb{G} |u(s,x)|^p   \mathrm{d}x \, \mathrm{d}s. 
\end{align*} If we introduce the time -- dependent functional
\begin{align*}
U_0(t)\doteq \int_{\mathbb{G}}  u(t,x) \, \mathrm{d}x,
\end{align*} then, we can rewrite the previous integral equality as follows:
\begin{align*}
U_0'(t)+U_0(t) - U_0'(0)-U_0(0)=\int_0^t\int_\mathbb{G} |u(s,x)|^p   \mathrm{d}x \, \mathrm{d}s \geqslant \int_0^t |U_0(s)|^p\, \mathrm{d}s,
\end{align*}  where in the last step we applied Jensen's inequality. 
 Multiplying the last inequality by $\mathrm{e}^{t}$, we have
\begin{align*}
\frac{\mathrm{d}}{\mathrm{d}t} (\mathrm{e}^{t} U_0(t)) = \mathrm{e}^{t} (U_0'(t)+U_0(t) )\geqslant (U_0'(0)+U_0(0))\, \mathrm{e}^{t}  +\mathrm{e}^{t} \int_0^t |U_0(s)|^p\, \mathrm{d}s.
\end{align*} Thus, integrating over $[0,t]$, we arrive at
\begin{align*}
\mathrm{e}^{t} U_0(t)\geqslant U_0(0) + (U_0'(0)+U_0(0))\, (\mathrm{e}^{t}  -1)+\int_0^t \mathrm{e}^{\tau} \int_0^\tau |U_0(s)|^p\, \mathrm{d}s \, \mathrm{d}\tau,
\end{align*} that is,
\begin{align} \label{fundamental ineq U0}
 U_0(t)\geqslant U_0(0)+ U_0'(0) (1-\mathrm{e}^{-t} )+\int_0^t \mathrm{e}^{\tau-t} \int_0^\tau |U_0(s)|^p\, \mathrm{d}s \, \mathrm{d}\tau.
\end{align} A first consequence of \eqref{fundamental ineq U0} is that $U_0$ is a positive function. Indeed, since $u_0,u_1$ are nonnegative and nontrivial functions, we have
\begin{align} \label{first lower bound U0}
U(t)\geqslant  U_0(0)+ U_0'(0) (1-\mathrm{e}^{-t} ) \geqslant C \varepsilon \qquad \mbox{for} \ t\geqslant 0,
\end{align} where the multiplicative constant $C$ depends on $u_0,u_1$. Furthermore, \eqref{fundamental ineq U0} provides the iteration frame
\begin{align} \label{iteration frame}
U(t)\geqslant  \int_0^t \mathrm{e}^{\tau-t} \int_0^\tau (U_0(s))^p\, \mathrm{d}s \, \mathrm{d}\tau.
\end{align}

\subsection{Iteration argument}
So far, we derived the iteration frame \eqref{iteration frame} and the first lower bound estimate for $U_0$ in \eqref{first lower bound U0}.
 Now we determine a sequence of lower bounds estimates for $U_0$ by using \eqref{iteration frame} in an iterative way. More precisely, we  show that
\begin{align}\label{sequence of lower bound U0}
U(t) \geqslant C_j\,  (t- L_j )^{\gamma_j} \qquad \mbox{for any} \ t\geqslant L_j,
\end{align} where $\{C_j\}_{j\in \mathbb{N}}$ and $\{\gamma_j\}_{j\in \mathbb{N}}$ are sequences of nonnegative real numbers that we be determined throughout this section and $\{L_j\}_{j\in\mathbb{N}}$ is the sequence of the partial products of the convergent infinite product
\begin{align*}
\prod_{k=0}^\infty \ell_k \quad  \mbox{with} \ \ \ell_k\doteq 1+p^{-k} \ \ \mbox{for any} \ k\in \mathbb{N},
\end{align*} that is, $$L_j\doteq \prod_{k=0}^{j} \ell_k \quad \mbox{for any} \ j\in \mathbb{N}.$$ 

Note that \eqref{first lower bound U0} implies \eqref{sequence of lower bound U0} for $j=0$ provided that $C_0\doteq C\varepsilon$ and $\gamma_0\doteq0$. Next we  prove \eqref{sequence of lower bound U0} by using an inductive argument. Consequently, it remains to show just the validity of the inductive step. Let us assume that \eqref{sequence of lower bound U0} is fulfilled for $j\geqslant 0$. We have to prove \eqref{sequence of lower bound U0} for $j+1$. After shrinking the domain of integration in \eqref{iteration frame}, if we plug \eqref{sequence of lower bound U0} in \eqref{iteration frame}, we get
\begin{align*}
U_0(t) & \geqslant \int_{L_j}^t \mathrm{e}^{\tau-t} \int_{L_j}^\tau (U_0(s))^p\, \mathrm{d}s \, \mathrm{d}\tau \\
& \geqslant C_j^p \int_{L_j}^t \mathrm{e}^{\tau-t} \int_{L_j}^\tau  (s- L_j )^{\gamma_j p}\, \mathrm{d}s \, \mathrm{d}\tau\\
 & \geqslant C_j^p (\gamma_jp+1)^{-1}  \int_{L_j}^t \mathrm{e}^{\tau-t} (\tau- L_j )^{\gamma_j p+1}\,  \mathrm{d}\tau\\
 & \geqslant C_j^p (\gamma_jp+1)^{-1}  \int_{t/\ell_{j+1}}^t \mathrm{e}^{\tau-t} (\tau- L_j )^{\gamma_j p+1}\,  \mathrm{d}\tau
\end{align*} for $t\geqslant L_{j+1} $. Note that in the last step we might restrict the domain of integration with respect to $\tau$ from $[L_j,t]$ to $[t/\ell_{j+1},t]$ since $t\geqslant L_{j+1}$ and $\ell_{j+1}>1$ imply $L_{j}  \leqslant t/\ell_{j+1}<t$. Therefore,
\begin{align*}
U_0(t) & \geqslant  C_j^p (\gamma_jp+1)^{-1} \ell_{j+1}^{-\gamma_jp-1}  (t- L_j \ell_{j+1} )^{\gamma_j p+1} \int_{t/\ell_{j+1}}^t \mathrm{e}^{\tau-t} \,  \mathrm{d}\tau \\
& \geqslant  C_j^p (\gamma_jp+1)^{-1} \ell_{j+1}^{-\gamma_jp-1}  (t- L_{j+1} )^{\gamma_j p+1}  \left(1-\mathrm{e}^{ -(1-1/\ell_{j+1})t} \right)
\end{align*} for $t\geqslant L_{j+1}$. Finally, we observe that for $t\geqslant L_{j+1}\geqslant \ell_{j+1} $ it holds the estimate 
\begin{align}
1-\mathrm{e}^{ -(1-1/\ell_{j+1})t} & \geqslant 1-\mathrm{e}^{-(\ell_{j+1}-1)} \geqslant 1- \left(1-(\ell_{j+1}-1)+\tfrac12 (\ell_{j+1}-1)^2\right) \notag \\
& = (\ell_{j+1}-1)\left(1-\tfrac12 (\ell_{j+1}-1)\right) 
 =  p^{-2(j+1)}\left(p^{j+1}-1/2\right)\geqslant \left(p-1/2\right)  p^{-2(j+1)}.
\end{align}
Hence, for $t\geqslant L_{j+1} $ we obtained 
\begin{align*}
U_0(t)  & \geqslant\left(p-1/2\right)   C_j^p (\gamma_jp+1)^{-1} \ell_{j+1}^{-\gamma_jp-1}  p^{-2(j+1)}  (t- L_{j+1} )^{\gamma_j p+1},
\end{align*}
which is precisely \eqref{sequence of lower bound U0} for $j+1$, provided that
\begin{align*}
C_{j+1} & \doteq   \frac{(p-1/2) \, C_j^p p^{-2(j+1)}}{(\gamma_j p+1) \,\ell_{j+1}^{\gamma_jp+1}}, \ \
\gamma_{j+1}  \doteq  1+p \gamma_j .
\end{align*}

\subsection{Upper bound estimate for the lifespan}

In the last section, we determined a sequence of lower bound estimates for $U$ in \eqref{sequence of lower bound U0}. Now we want to show that the $j$ -- dependent lower bound in \eqref{sequence of lower bound U0} for $U_0$ blows up  as $j\to \infty$ for $t$ greater than a certain $\varepsilon$ -- dependent threshold. This will prove the desired blow~--~up result and, as byproduct, will provide the upper bound estimate for the lifespan. Let us get started by estimating the multiplicative constant $C_j$ from below.

In order to estimate $C_j$ we need to determine first the explicit representation for $\gamma_{j}$.

 Since  $\gamma_j= 1+p\gamma_{j-1}$, applying recursively this relation, we get
\begin{align}
\gamma_j &= p^2\gamma_{j-2}+1+p = \cdots =   p^{j}\gamma_0+1+p+\cdots +p^{j-1}=\tfrac{p^{j}-1}{p-1}. \label{representation gamma j}
\end{align} Thus, from \eqref{representation gamma j} it results
\begin{align*}
\gamma_{j-1}p+1 = \gamma_j \leqslant \tfrac{p^{j}}{p-1}
\end{align*} which yields in turns
\begin{align*}
C_j \geqslant  (p-1/2) (p-1)C_{j-1}^p p^{-3j} \ell_{j}^{-\gamma_{j}}.
\end{align*}
Moreover, it holds
\begin{align*}
\lim_{j\to \infty} \ell_j^{\gamma_j}= \lim_{j\to \infty} \exp\left(\frac{p^{j}}{p-1}\log \left(1+p^{-j}\right)  \right) = \mathrm{e}^{1/(p-1)}.
\end{align*} In particular, there exists a suitable constant $M=M(p)>0$ such that $\ell_j^{-\gamma_j}\geqslant M$ for any $j\in \mathbb{N}$. Hence, 
\begin{align*}
C_j \geqslant \underbrace{ (p-1/2)(p-1)  M }_{\doteq K}C_{j-1}^p p^{-3j} \quad \mbox{for any} \ j\in\mathbb{N}.
\end{align*}
Applying the logarithmic function to both sides of the inequality $C_j\geqslant K\, C_{j-1}^p p^{-3j}$ and using iteratively this inequality, we find
\begin{align*}
\log C_j & \geqslant p \log C_{j-1} -3 j \log p +\log K \\
& \geqslant p^2 \log C_{j-2} -3( j+(j-1)p ) \log p + (1+p)\log K \\
& \geqslant \cdots \geqslant p^{j}\log C_0 -3 \left(\sum_{k=0}^{j-1} (j-k)p^k\right) \log p+ \left(\sum_{k=0}^{j-1} p^k\right)\log K.
\end{align*} Employing the identities
\begin{align} \label{identity sum (j-k)p^k}
\sum_{k=0}^{j-1} (j-k)p^k = \frac{1}{p-1} \left(\frac{p^{j+1}-p}{p-1}-j\right)\quad\text{and}\quad  \sum_{k=0}^{j-1} p^k = \frac{p^j-1}{p-1},
\end{align} it follows
\begin{align*}
\log C_j \geqslant p^j \left(\log C_0-\frac{3p \log p}{(p-1)^2} +\frac{\log K}{p-1} \right)+\frac{3j \log p}{p-1}+\frac{3p \log p}{(p-1)^2}-\frac{\log K}{p-1}
\end{align*} for any $j\in\mathbb{N}$. Let $j_0=j_0(p)\in\mathbb{N}$ be the smallest nonnegative integer such that
\begin{align*}
j_0\geqslant \frac{\log K}{3\log p}-\frac{p}{p-1}.
\end{align*}
Therefore, for any $j\geqslant j_0$ we have
\begin{align}
\log C_j &\geqslant p^j \left(\log C_0-\frac{3p \log p}{(p-1)^2} +\frac{\log K}{p-1} \right) = p^j \log \left(K^{1/(p-1)}p^{-(3p)/(p-1)^2}C_0\right) \notag \\ & =p^j \log (E_0 \varepsilon) \label{lower bound log Cj}
\end{align} where $E_0=K^{1/(p-1)}p^{-(3p)/(p-1)^2}C>0$.

Let us denote $$L\doteq \lim_{j\to\infty} L_j = \prod_{j=0}^\infty \ell_{j}\in \mathbb{R}.$$ Notice that due to $\ell_j>1$, it results $L_j \uparrow L$ as $j\to \infty$. In particular, \eqref{sequence of lower bound U0} holds for any $j\in\mathbb{N}$ and any $t\geqslant L$.

Combining \eqref{sequence of lower bound U0}, \eqref{representation gamma j} and \eqref{lower bound log Cj}, we find
\begin{align*}
U_0(t) & \geqslant \exp \left(p^j \log (E_0\varepsilon)\right)  (t- L)^{\gamma_j} \\
& = \exp \left(p^j  \left( \log (E_0\varepsilon)+\tfrac{1}{p-1}\log(t-L)\right)\,\right) (t- L)^{-1/(p-1)}
\end{align*}
for any $j\geqslant j_0$ and any $t\geqslant L$. Finally, for $t\geqslant 2L$, since $t-L \geqslant t/2$, we have
\begin{align} \label{final lower bound U0}
U_0(t) & \geqslant  \exp \left(p^j \log \left(E_1\varepsilon \, t^{\frac{1}{p-1}}\right)\right)  (t- L)^{-1/(p-1)}
\end{align}
for any $j\geqslant j_0$, where $E_1\doteq 2^{-1/(p-1)}E_0$. For any $p>1$, the exponent for $t$ in the logarithmic term in \eqref{final lower bound U0} is positive.
 Let us fix $\varepsilon_0~=~\varepsilon_0(u_0,u_1,p)>0$ such that
\begin{align*}
\varepsilon_0\leqslant (2L)^{-1/(p-1)} E_1^{-1}.
\end{align*} Consequently, for any $\varepsilon\in(0,\varepsilon_0]$ and any $t>(E_1 \varepsilon)^{-(p-1)}$, we have
\begin{align*}
t\geqslant 2L \quad \mbox{and} \quad \log \left(E_1\varepsilon \, t^{\frac{1}{p-1}}\right)>0.
\end{align*} Thus, for any $\varepsilon\in(0,\varepsilon_0]$ and any $t>(E_1 \varepsilon)^{-(p-1)}$ letting $j\to \infty$ in \eqref{final lower bound U0} we see that the lower bound for $U_0(t)$ blows up. So, for any $\varepsilon\in(0,\varepsilon_0]$ the functional $U_0$ has to blow up in finite time and, moreover, the lifespan of the local solution $u$ can be estimated from above as follows: $$T(\varepsilon) \lesssim \varepsilon^{-(p-1)}.$$
We completed the proof of Theorem \ref{Thm blow up}.

\section{Final remarks}

In this paper we considered the Laplace -- Beltrami operator on $\mathbb{G}$. Nonetheless, we might study the damped wave operator $\partial_t^2-\mathcal{L}_X+\partial_t$ with a sub -- Laplacian $\mathcal{L}_X$ associated to a system $X=\{X_1,\ldots,X_k\}$ of left -- invariant vector fields satisfying H\"ormander condition in place of the Laplace -- Beltrami operator. 

In this case, the proof of the blow -- up result is exactly the same, as the sub -- Laplacian is formally self -- adjoint, and we can introduce energy solutions analogously as in Definition \ref{Definition energy sol}.

 For the proof a local in time result we may proceed similarly as in Section \ref{Section local existence} keeping in mind that $\sigma_{\mathcal{L}_X}(\xi)$ can be written in diagonal form for any $[\xi]\in\widehat{\mathbb{G}}$ (cf. \cite[Formula (3.3)]{GR15}) and that the upper bound for $p$ due to a Gagliardo -- Nirenberg type inequality is $d/(d-2)$ for $d\geqslant 3$, where $d=d(X)$ is the local dimension of $\mathbb{G}$ with respect to the system $X$ (see also Remark \ref{Remark local dimension}).

As we pointed out in the introduction, we might interpret the results obtained in this paper by saying that the critical exponent is the Fujita exponent in the 0 -- dimensional case. So, rather than the topological dimension of the group $\mathbb{G}$, what determines the critical exponent is the global dimension of $\mathbb{G}$, which is 0 for a compact Lie group.

In a series of forthcoming papers \cite{Pal20WE,Pal20WEdm}, other semilinear hyperbolic models will be considered on compact Lie groups. 



\end{document}